\documentclass{article}
\usepackage{geometry}
\usepackage{IEEEtrantools}
\usepackage{amsmath}
\usepackage{amssymb}
\usepackage{bm}
\usepackage{acronym}
\usepackage{authblk}
\usepackage{hyperref}

\acrodef{DMD}{Dynamic Mode Decomposition}
\acrodef{SVD}{Singular Value Decomposition}
\acrodef{POD}{Proper Orthogonal Decomposition}

\begin{document}

\title{Dynamic Mode Decomposition and Koopman Theory}
\author{Sourya Dey}
\affil{Galois, Inc.}
\date{}

\maketitle

\section*{Abstract}
\acf{DMD} is a technique to approximate generally non-linear dynamical systems using linear techniques, which are better understood and easier to analyze. Koopman theory extends \ac{DMD} by transforming the original system into a new domain which facilitates linearization. This is a technical report on \ac{DMD} and Koopman theory, with primary focus on explaining the underlying mathematics in clear and concise form. We include dimensions of vectors and marices, and step-by-step derivations of equations in order to assist the user in easily comprehending these concepts. This report will also enable users to implement \ac{DMD} and Koopman theory in code.

\section*{Acknowledgements}
The material in this report has been compiled from relevant portions of \cite{dmdbook,Tu2014}.

\section*{Notation}
Matrices are written in bold capital letters, e.g. $\bm{A}$, and vectors in bold small letters, e.g. $\bm{a}$. Dimensions of quantities are written under them, e.g. $\underset{m\times n}{\bm{A}}$ denotes a matrix $\bm{A}$ with $m$ rows and $n$ columns.

\section{Dynamic Mode Decomposition}\label{sec-dmd}
Consider a system whose state at index $k$ is described by a vector $\bm{x}_k$. Essentially $\bm{x}_k$ is a vector of observed values sampled at index $k$. The index can be time, or any other independent variable. Similarly, the sampled state at index $k+1$ is $\bm{x}_{k+1}$. $\bm{x}$ has a finite number of dimensions equal to the number of observed values (assume $n$). This means that it can be expressed as a linear combination of $n$ basis vectors.

The transition function to take $\bm{x}_k$ to $\bm{x}_{k+1}$ is $F$, i.e.:
\begin{IEEEeqnarray}{c}\label{eq-original-nonlinear-vector-transition}
\bm{x}_{k+1} = F(\bm{x}_k)
\end{IEEEeqnarray}

$F(\cdot)$ is, in general, a non-linear function, i.e. $F(a\bm{x}_1+b\bm{x}_2) \ne aF(\bm{x}_1)+bF(\bm{x}_2)$.\footnote{The function $F(\cdot)$ is written using a capital letter because typically the space is described as $\frac{d\bm{x}}{dt} = f(\bm{x})$. Integrating this from $t=k$ to $t=k+1$ gives $\bm{x}_{k+1} = F(\bm{x}_k)$.} \textbf{Our goal is to find a matrix $\bm{A}$ such that the transition dynamics become linear}, i.e.:
\begin{IEEEeqnarray}{c}\label{eq-original-linear-vector-transition}
\underset{n\times1}{\bm{x}_{k+1}} \approx \underset{n\times n}{\bm{A}}\underset{n\times1}{\bm{x}_k}
\end{IEEEeqnarray}

This is extremely powerful. If this indeed holds, then we can use the eigendecomposition of $\bm{A}$ to get the state $\bm{x}_j$ at any index $j$. Assume the eigenvectors of $\bm{A}$ are $\bm{W}$, and the eigenvalues are $\bm{\Lambda}$, i.e. $\underset{n\times n}{\bm{A}}\underset{n\times n}{\bm{W}} = \underset{n\times n}{\bm{W}}\underset{n\times n}{\bm{\Lambda}}$. These eigenvectors are called the \ac{DMD} modes and they form a basis in $n$-dimensional space, such that any vector $x$ can be expressed as $\bm{x} = \sum_{i=1}^nb_i\bm{w}_i$, where the $b$ values are the coefficients of the linear combination. Putting them together in a vector $\bm{b}$, we get $\bm{x} = \bm{W}\bm{b}$.

Assume the system has known initial state $\bm{x}_0 = \bm{W}\bm{b}$. We can solve for $\bm{b}$ as:
\begin{IEEEeqnarray}{c}
\underset{n\times1}{\bm{b}} = \underset{n\times n}{\bm{W}^{-1}}\underset{n\times1}{\bm{x}_0}
\end{IEEEeqnarray}

Then, we get $\bm{x}_1 = \bm{A}\bm{W}\bm{b} = \bm{W}\bm{\Lambda}\bm{b}$. Extending this to any index, we get the incredibly powerful:
\begin{IEEEeqnarray}{c}
\underset{n\times1}{\bm{x}_k} = \underset{n\times n}{\bm{W}}\underset{n\times n}{\bm{\Lambda}^k}\underset{n\times1}{\bm{b}} = \sum_{i=1}^n{\lambda_i}^kb_i\underset{n\times1}{\bm{w}_i}
\end{IEEEeqnarray}
\textbf{for any index $k$}, including positive and negative values of high magnitude. When the index is time, this implies that we can compute states of the system way out into the future or back in the past. Note, however, that this means that \textbf{for a stable system, the eigenvalue with largest real part (i.e. dominant eigenvalue) should have real part equal to $1$}. If it exceeds $1$, the system will explode in the future and implode in the past. If it's less than $1$, the system will explode in the past and implode in the future.

Note that ideally \emph{all eigenvalues} (not just the dominant one) should lie on the unit circle, because they all affect the system (albeit to a lesser extent than the dominant eigenvalue). If the non-dominant eigenvalues are close to the origin (i.e. real part much less than $1$), then the system still has a chance of exploding in the past. Such numerical instability should be kept in mind when predicting states of the system in the past.

\section{Implementing DMD}
More generally, suppose there exist $m+1$ measurements of $\bm{x}$ at different indexes, i.e. $\bm{x}_0$ to $\bm{x}_m$. These are collected as $n\times m$ matrices $\bm{X} = \left[\bm{x}_0, \cdots, \bm{x}_{m-1}\right]$ and $\bm{X}' = \left[\bm{x}_1, \cdots, \bm{x}_m\right]$. Then, we get the matrix form of Eq. \eqref{eq-original-linear-vector-transition}:
\begin{IEEEeqnarray}{c}\label{eq-original-linear-matrix-transition}
\underset{n\times m}{\bm{X}'} \approx \underset{n\times n}{\bm{A}}\underset{n\times m}{\bm{X}}
\end{IEEEeqnarray}

There can be 2 different cases:
\subsection*{Case 1: $m\geq n$}\label{sec-dmd-mbig}
This is the standard regression case where we have more observations than features. We cannot solve for $\bm{A}$ exactly since we have $mn$ equations and only $n^2$ variables. So, we try to minimize the least squares error $\left\|\bm{X}'-\bm{A}\bm{X}\right\|^2$, which gives the solution:
\begin{IEEEeqnarray}{c}
\underset{n\times n}{\bm{A}} = \underset{n\times m}{\bm{X}'}\underset{m\times n}{\bm{X}^{\dagger}}
\end{IEEEeqnarray}

\subsection*{Case 2: $m<n$}\label{sec-dmd-nbig}
This is an over-determined system where we are trying to solve for $n^2$ variables in only $mn$ equations. There can be infinite solutions. Whatever solution for $A$ we pick, it will be an $n\times n$ matrix with rank only $m$.

\paragraph{\acf{SVD}:} Consider the \ac{SVD} of $\bm{X}$:
\begin{IEEEeqnarray}{c}
\underset{n\times m}{\bm{X}} = \underset{n\times n}{\bm{U}}\underset{n\times m}{\bm{\Sigma}}\underset{m\times m}{\bm{V}^H}
\end{IEEEeqnarray}
where $^H$ is the conjugate transpose or Hermitian.

\paragraph{Truncated \ac{SVD}:} Since $m<n$, the rank of $\bm{X}$ is at most $m$. Let's assume the rank of $\bm{X}$ is actually $r$, which is less than $m$. (If it turns out that the rank of $\bm{X}$ is actually $m$, replace $r$ with $m$ in the following). Now, consider the truncated \ac{SVD}. Since the rank is $r$, we need only keep the left $r$ columns of both $\bm{U}$ and $\bm{V}$, and the top-left $r\times r$ portion of $\bm{\Sigma}$. Then we get:
\begin{IEEEeqnarray}{c}
\underset{n\times m}{\bm{X}} = \underset{n\times r}{\bm{U}}\underset{r\times r}{\bm{\Sigma}}\underset{r\times m}{\bm{V}^H}
\end{IEEEeqnarray}
The left singular vectors $\bm{U}$ are the \ac{POD} modes.

\paragraph{Reduced dimensionality \ac{POD} space:} The original $n$-dimensional vectors $\bm{x}$ can be projected into the reduced $r$-dimensional space of \ac{POD} modes using $\tilde{\bm{x}} = \bm{U}^H\bm{x}$. Doing this to the whole matrix at once gives:
\begin{IEEEeqnarray}{c}
\underset{r\times m}{\tilde{\bm{X}}} = \underset{r\times n}{\bm{U}^H}\underset{n\times m}{\bm{X}}
\end{IEEEeqnarray}
and to convert back, we use $\bm{x} = \bm{U}\tilde{\bm{x}}$ and $\bm{X} = \bm{U}\tilde{\bm{X}}$.

Since $r$ has replaced $n$ and $m\geq r$, we are back in Case 1. We will now proceed with this case.

\subsection{The DMD algorithm}\label{sec-dmd-algo}
\paragraph{Finding $\bm{A}$:} Using the \ac{SVD} of $\bm{X}$, we can get $\bm{A} = \bm{X}'\bm{X}^{\dagger} = \bm{X}'\bm{V}\bm{\Sigma}^{-1}\bm{U}^H$.

We can also get $\tilde{\bm{A}} = \tilde{\bm{X}}'\tilde{\bm{X}}^{\dagger}$. This is where it gets tricky. $\tilde{\bm{X}}'$ is the projection of $\bm{X}'$ on the same \ac{POD} modes obtained from the \ac{SVD} of $\bm{X}$, i.e. $\tilde{\bm{X}}' = \bm{U}^H\bm{X}'$. So we get:
\begin{IEEEeqnarray}{rCl}
\underset{r\times r}{\tilde{\bm{A}}} &=& \underset{r\times m}{\tilde{\bm{X}}'}\underset{m\times r}{\tilde{\bm{X}}^{\dagger}} \IEEEyesnumber\IEEEyessubnumber\\
&=& \underset{r\times n}{\bm{U}^H}\underset{n\times m}{\bm{X}'}\underset{m\times n}{\bm{X}^{\dagger}}\underset{n\times r}{\bm{U}} \IEEEnonumber\IEEEyessubnumber\\
&=& \underset{r\times n}{\bm{U}^H}\underset{n\times n}{\bm{A}}\underset{n\times r}{\bm{U}} \IEEEnonumber\IEEEyessubnumber\\
&=& \underset{r\times n}{\bm{U}^H}\underset{n\times m}{\bm{X}'}\underset{m\times r}{\bm{V}}\underset{r\times r}{\Sigma^{-1}} \IEEEnonumber\IEEEyessubnumber
\end{IEEEeqnarray}
\textbf{This means that we never have to compute $\bm{A}$, we can just work with the reduced dimensionality $\tilde{\bm{A}}$. The reduced dimension $r$ can be chosen to be anything as long as it's upper-bounded by $\min(n,m)$.}

\paragraph{Eigendecomposition} We can now compute the eigendecomposition:
\begin{IEEEeqnarray}{c}
\underset{r\times r}{\tilde{\bm{A}}} = \underset{r\times r}{\tilde{\bm{W}}}\underset{r\times r}{\tilde{\bm{\Lambda}}}\underset{r\times r}{\tilde{\bm{W}}^{-1}}
\end{IEEEeqnarray}
The first $r$ eigenvalues of $\bm{A}$ are the same as those of $\tilde{\bm{A}}$, i.e. $\bm{\Lambda} = \tilde{\bm{\Lambda}}$. Since the rank of $\bm{X}$ and $\bm{X}'$ is $r$, the rank of $\bm{A}$ is also $r$ (since multiplying two rank $r$ matrices makes the result rank $r$). So, the other $n-r$ eigenvalues of $\bm{A}$ are all $0$.

This means that we only care about the first $r$ eigenvectors of $\bm{A}$. The eigenvectors $\bm{W}$ of $\bm{A}$ can be found in 2 ways:

The \emph{exact} eigenvectors are:
\begin{IEEEeqnarray}{c}
\underset{n\times r}{\bm{W}} = \underset{n\times m}{\bm{X}'}\underset{m\times r}{\bm{V}}\underset{r\times r}{\Sigma^{-1}}\underset{r\times r}{\tilde{\bm{W}}}
\end{IEEEeqnarray}
For a backward proof, we can compute $\bm{W}\bm{\Lambda}\bm{W}^{\dagger}$ to get back $\bm{A}$.

The \emph{projected} eigenvectors are simply:
\begin{IEEEeqnarray}{c}
\underset{n\times r}{\bm{W}} = \underset{n\times r}{\bm{U}}\underset{r\times r}{\tilde{\bm{W}}}
\end{IEEEeqnarray}
which is what one would expect given the rules of the \ac{POD} space are such that $\bm{X} = \bm{U}\tilde{\bm{X}}$. Note that this $\bm{W}$ is equal to $\bm{X}\bm{V}\bm{\Sigma}^{-1}\tilde{\bm{W}}$, which is equal to the exact formulation except that $\bm{X}'$ is replaced with $\bm{X}$. Thus, the projected formulation will be equal to the exact formulation if $\bm{X}$ and $\bm{X}'$ have the same column spaces.

Note that computing the projected eigenvectors is more numerically stable. However, it is recommended to use the exact formulation unless a particular eigenvalue and its corresponding exact eigenvector are both $0$, in which case the projected eigenvector should be used.

\paragraph{Calculating other $\bm{x}$ values:} The rest of the process is mostly the same. We find:
\begin{IEEEeqnarray}{c}\label{eq-b}
\underset{r\times1}{\bm{b}} = \underset{r\times n}{\bm{W}^{\dagger}}\underset{n\times1}{\bm{x}_0}
\end{IEEEeqnarray}

Then we can find states for any index $k$ as:
\begin{IEEEeqnarray}{c}\label{eq-discrete-solution}
\underset{n\times1}{\bm{x}_k} = \underset{n\times r}{\bm{W}}\underset{r\times r}{\bm{\Lambda}^k}\underset{r\times1}{\bm{b}} = \sum_{i=1}^r{\lambda_i}^kb_i\underset{n\times1}{\bm{w}_i}
\end{IEEEeqnarray}

\subsection{Converting between discrete and continuous indexes}
Thus far we have dealt with discrete indexes, i.e. the vector $\bm{x}$ was sampled at indexes with integral $k$. What if we want to find a continuous representation $\bm{x}(k)$, where $k$ can be non-integral? For example, we may wish to find the state $\bm{x}$ of the system at index $3.5$. In this case, the dynamics of the system are described as:
\begin{IEEEeqnarray}{c}
\underset{n\times1}{\frac{d\bm{x}}{dk}} = \underset{n\times n}{\bm{\mathcal{A}}}\underset{n\times1}{\bm{x}}
\end{IEEEeqnarray}

This has solution:
\begin{IEEEeqnarray}{c}
\bm{x}(k) = e^{\bm{\mathcal{A}}k}\bm{x}_0
\end{IEEEeqnarray}

The eigendecomposition is $\bm{\mathcal{A}} = \bm{W}\bm{\Omega}\bm{W}^{-1}$. Recall that any function of $\bm{\mathcal{A}}$ gets applied to the eigenvalues, leaving the eigenvectors intact. As before, $\bm{x}_0 = \bm{W}\bm{b}$.

Then we can solve for any index as:
\begin{IEEEeqnarray}{c}
\underset{n\times1}{\bm{x}(k)} = \bm{W}e^{\bm{\Omega}k}\bm{W}^{-1}\bm{x}_0 = \underset{n\times n}{\bm{W}}\underset{n\times n}{e^{\bm{\Omega}k}}\underset{n\times1}{\bm{b}}
\end{IEEEeqnarray}
or, for the truncated case:
\begin{IEEEeqnarray}{c}\label{eq-continuous-solution}
\underset{n\times1}{\tilde{\bm{x}}(k)} = \underset{n\times r}{\bm{W}}\underset{r\times r}{e^{\bm{\Omega}k}}\underset{r\times1}{\bm{b}}
\end{IEEEeqnarray}

This formulation is equivalent to the discrete index case if we write $\bm{A} = e^{\bm{\mathcal{A}}\Delta k}$, where $\Delta k$ is the sampling interval. \textbf{This means that the eigenvectors of both matrices are the same $\bm{W}$, and the eigenvalues are related as:}
\begin{IEEEeqnarray}{rCl}\label{eq-eigenvalue-conversion}
\lambda &=& e^{\omega\Delta k} \IEEEyesnumber\IEEEyessubnumber\\
\Rightarrow \omega &=& \frac{\ln{\lambda}}{\Delta k} \IEEEnonumber\IEEEyessubnumber
\end{IEEEeqnarray}

As an example, suppose we want to find $\bm{x}$ at index $3.5$. For the continuous index case, we would get $\bm{x}(3.5) = \bm{W}e^{3.5\bm{\Omega}}\bm{b}$. For the usual discrete sampling case with $\Delta k=1$, we cannot find $\bm{x}_{3.5}$ since the index $3.5$ isn't a valid sampling index. However, we can solve for $\bm{x}$ at index $3.5$ if we take $\Delta k = 0.5$. This means that we need to retake the measurements so that $\bm{X}'$, which stands for one discrete sampling index ahead of $\bm{X}$, will now have an actual value that is $0.5$ ahead. Then, $\bm{\Lambda} = e^{0.5\bm{\Omega}}$, so $e^{3.5\bm{\Omega}} = \bm{\Lambda}^7$. This makes sense since $3.5$ is the $7$th index in this new sampling scheme. Then we can solve $\bm{x}(3.5) = \bm{x}_7 = \bm{W}\bm{\Lambda}^7\bm{b}$.

\paragraph{Tying it together:} Discrete indexing is obviously the only way we can operate on a computer. So, assume we have state measurements at indexes $[0,0.5,1,1.5,2,3,4,5]$, and we want to predict the states at indexes $2.5$ and $8$. We can do it in one out of two ways:
\begin{itemize}
    \item Construct the discrete problem with $\Delta k = 0.5$. This is not possible for the complete data since the values at $2.5$, $3.5$ and $4.5$ are unavailable. It becomes possible if we only calculate \ac{DMD} using indexes $[0,0.5,1,1.5,2]$, which, in this new scheme, become indexes $[0,1,2,3,4]$. Construct $\bm{X} = [\bm{x}_0\cdots\bm{x}_3]$ (which is actually $[\bm{x}(0)\cdots\bm{x}(1.5)]$), and $\bm{X}' = [\bm{x}_1\cdots\bm{x}_4]$ (which is actually $[\bm{x}(0.5)\cdots\bm{x}(2)]$), and proceed as before. Once the eigendecomposition of the system's $\bm{A}$ has been determined, calculate $\bm{x}(2.5)$ as $\bm{x}_5$ and $\bm{x}(8)$ as $\bm{x}_{16}$ using Eq. \eqref{eq-discrete-solution}.
    \item Construct the discrete problem with $\Delta k = 1$. Only consider indexes $[0,1,2,3,4,5]$ for calculating \ac{DMD}. Construct $\bm{X} = [\bm{x}_0\cdots\bm{x}_4]$, $\bm{X}' = [\bm{x}_1\cdots\bm{x}_5]$, and proceed as before. Once the eigendecomposition of the system's $\bm{A}$ has been determined, convert the discrete eigenvalues $\bm{\Lambda}$ to continuous eigenvalues $\bm{\Omega}$ using Eq. \eqref{eq-eigenvalue-conversion}, and calculate $\bm{x}(2.5)$ and $\bm{x}(5)$ using Eq. \eqref{eq-continuous-solution}.
\end{itemize}

\subsection{Different starting index}
Thus far we have assumed that the starting index is $\bm{x}_0$. What if it's not? Suppose we have measurements starting at $\bm{x}_i$, $i\ne0$. There are 2 methods to deal with this:
 
\paragraph{Assume index $0$ exists:} Assume there is an un-measured $\bm{x}_0$. Then, $\bm{x}_i = \bm{A}^i\bm{x}_0$. Putting the earlier expression $\bm{x}_0 = \bm{W}\bm{b}$ into this yields:
\begin{IEEEeqnarray}{rCl}
\bm{x}_i &=& \left(\bm{W}\bm{\Lambda}\bm{W}^{-1}\right)^i \bm{W}\bm{b} \IEEEnonumber \\
&=& \bm{W}\bm{\Lambda}^i\bm{b} \IEEEnonumber \\
\Rightarrow \bm{b} &=& \bm{\Lambda}^{-i}\bm{W}^{-1}\bm{x}_i \IEEEyesnumber
\end{IEEEeqnarray}
In terms of the specific notation used so far, Eq. \eqref{eq-b} now becomes:
\begin{IEEEeqnarray}{c}
\underset{r\times1}{\bm{b}} = \underset{r\times r}{\bm{\Lambda}^{-i}}\underset{r\times n}{\bm{W}^{\dagger}}\underset{n\times1}{\bm{x}_i}
\end{IEEEeqnarray}
Thus, we can get $\bm{b}$ in terms of whichever index $\bm{x}$ starts from by using the appropriate opposite power on the eigenvalue matrix. We don't need to know $\bm{x}_0$.

Then, any other value $\bm{x}_j$ can be obtained as:
\begin{IEEEeqnarray}{rCl}
\bm{x}_j &=& \bm{W}\bm{\Lambda}^j\bm{b} \IEEEnonumber \\
&=& \bm{W}\bm{\Lambda}^j\bm{\Lambda}^{-i}\bm{W}^{\dagger}\bm{x}_i \IEEEnonumber \\
&=& \bm{A}^{j-i}\bm{x}_i \IEEEyesnumber
\end{IEEEeqnarray}
which is as one would expect when computing index $j$ from index $i$. Also note that $j$ can be less than $i$, and everything will still hold. For example, $\bm{x}_0 = \bm{A}^{-i}\bm{x}_i$.

\paragraph{Shift the indexes so that $i$ becomes $0$:} Thus just means subtracting $i$ from any index. So, the given $\bm{x}_i$ in the original index space becomes $\bm{x}_{i-i} = \bm{x}_0$ in the shifted index space, and $\bm{x}_j$ becomes $\bm{x}_{j-i}$. That way, we retain the relation $\bm{b} = \bm{W}^{-1}\bm{x}_0$ in the shifted index space, i.e. we avoid powers on the eigenvalue matrix when computing $\bm{b}$ and can keep the math the same.

Now, suppose we want to compute $\bm{x}_j$ in the original index space. Then, in the shifted index space, we need to compute:
\begin{IEEEeqnarray}{c}
\bm{x}_{j-i} = \bm{W}\bm{\Lambda}^{j-i}\bm{b}
\end{IEEEeqnarray}
as one would expect. Also note that $\bm{x}_0$ in the original index space will be computed in the shifted index space as:
\begin{IEEEeqnarray}{rCl}
\bm{x}_{0-i} &=& \bm{W}\bm{\Lambda}^{0-i}\bm{b} \IEEEnonumber \\
&=& \bm{W}\bm{\Lambda}^{0-i}\bm{W}^{-1}\bm{x}_0 \IEEEnonumber \\
&=& \bm{A}^{-i}\bm{x}_0 \IEEEnonumber \\
&=& \bm{A}^{-i}\bm{x}_{i-i} \IEEEnonumber
\end{IEEEeqnarray}
When expressed in the original index space, this becomes $\bm{x}_0 = \bm{A}^{-i}\bm{x}_i$, which is identical to the previous method. Thus, the two methods give identical results, as they should.

\subsection{Multiple trajectories}
A trajectory is defined as a run of a dynamical system from some initial state to some final state, i.e. $\left[\bm{x}_0, \cdots, \bm{x}_m\right]$. So far, we have been dealing with single trajectories. But it can so happen that a dynamical system is described by multiple trajectories instead of a single trajectory. For example, consider the system $\bm{x}_{k+1} = \bm{x}_k+\psi$, where $\psi$ denotes noise. We want a single $\bm{A}$ for the whole system since its dynamics don't change. But we also need a way to represent data from the multiple trajetories. This can be done by \emph{concatenating the trajectories along the index dimension}.

Suppose we perform $T$ runs of a system, which have $\left[m_1+1, m_2+1, \cdots, m_T+1\right]$ indexed states, respectively. Using superscript for trajectory, the states can be described as $\left[\bm{x}^1_0, \cdots, \bm{x}^1_{m_1}\right], \left[\bm{x}^2_0, \cdots, \bm{x}^2_{m_2}\right], \cdots, \left[\bm{x}^T_0, \cdots, \bm{x}^T_{m_T}\right]$. Concatenating along the index dimension, we get:
\begin{IEEEeqnarray}{rCl}
\bm{X} &=& \left[\bm{x}^1_0, \cdots, \bm{x}^1_{m_1-1}, \quad \bm{x}^2_0, \cdots, \bm{x}^2_{m_2-1}, \quad \cdots, \quad \bm{x}^T_0, \cdots, \bm{x}^T_{m_T-1}\right] \IEEEyesnumber \IEEEyessubnumber \\
\bm{X}' &=& \left[\bm{x}^1_1, \cdots, \bm{x}^1_{m_1}, \quad \bm{x}^2_1, \cdots, \bm{x}^2_{m_2}, \quad \cdots, \quad \bm{x}^T_1, \cdots, \bm{x}^T_{m_T}\right] \IEEEnonumber \IEEEyessubnumber
\end{IEEEeqnarray}
Both these matrices have dimensions $n\times M$, where $M = \sum_{t=1}^{T}{m_t}$. $\bm{X}$ leaves out the last sample of each trajectory, while $\bm{X}'$ leaves out the first sample of each trajectory.

The rest of the algorithm proceeds as given previously, except that $m$ is replaced by $M$. This means that the reduced dimension (i.e. rank) $r$ will be upper-bounded by $\min(n,M)$.

\section{Koopman theory}
The problem with \ac{DMD} is that the original dynamics may not be linearizable. However, what if we transform the problem from the existing $n$-dimensional space to a different space? We can replace all vectors $\bm{x}$ with functions $g(\bm{x})$. Thus, $g(\bm{x})$ is an \textbf{encoding}. Essentially, now the encoded value(s) at index $k$ is $g(\bm{x}_k)$ and at index $k+1$ is $g(\bm{x}_{k+1})$. Since $g(\bm{x})$  is a function, in general, it has an infinite number of dimensions, i.e. the number of basis functions whose linear combination equals $g(\bm{x})$  is infinite.

The transition operator to take g($\bm{x}_k$) to g($\bm{x}_{k+1}$) is $\bm{K}$, which is the Koopman operator. So:
\begin{IEEEeqnarray}{c}\label{eq-koopman-vector-transition}
g(\bm{x}_{k+1}) = \bm{K}g(\bm{x}_k)
\end{IEEEeqnarray}

We assume that $\bm{K}$ is a linear operator, i.e. it is a (infinite-dimensional) matrix. \textbf{Essentially, the Koopman operator transforms the problem of a non-linear transition function acting in a finite-dimensional vector space to the problem of a linear operator acting on an infinite-dimensional function space.}

The basis functions for this encoded space can be the eigenfunctions of $\bm{K}$, i.e. the functions $\{\phi(\cdot)\}$ which satisfy $\bm{K}\phi(\bm{x}) = \lambda\phi(\bm{x})$. So, any function $g(\cdot)$ on which $\bm{K}$ operates can be written as a linear, infinite combination of the Koopman eigenfunctions $\{\phi(\cdot)\}$.

\subsection{Koopman operator as function composition}
Note that $\bm{K}g(\bm{x}_k) = g(\bm{x}_{k+1})$. But as per Eq. \eqref{eq-original-nonlinear-vector-transition}, $\bm{x}_{k+1} = F(\bm{x}_k)$. This implies that $\bm{K}g(\bm{x}_k) = g(F(\bm{x}_k))$. Dropping the index:
\begin{IEEEeqnarray}{c}
\bm{K}g(\bm{x}) = g(F(\bm{x}))
\end{IEEEeqnarray}
So, the Koopman operator can be regarded as function composition.

\section{Implementing Koopman theory as generalized DMD}
The Koopman operator cannot be implemented in practice since it is infinite dimensional. However, we can implement a finite approximation by assuming the encoded space to be $p$-dimensional, where $p$ is finite. Then, we can write:
\begin{IEEEeqnarray}{c}
\underset{p\times1}{\bm{y}} = g\left(\underset{n\times1}{\bm{x}}\right)
\end{IEEEeqnarray}
and use a finite-dimensional $p\times p$ Koopman matrix $\bm{K}$ as the equivalent of the $n\times n$ \ac{DMD} matrix $\bm{A}$ from Eq. \eqref{eq-original-linear-vector-transition}:
\begin{IEEEeqnarray}{c}
\underset{p\times1}{\bm{y}_{k+1}} \approx \underset{p\times p}{\bm{K}}\underset{p\times1}{\bm{y}_k}
\end{IEEEeqnarray}

Thus, once we have the $\bm{y}$ values, \textbf{one can perform \ac{DMD} exactly as described thus far by replacing $\bm{x}$ with $\bm{y}$, $\bm{X}$ with $\bm{Y}$, $n$ with $p$, $\bm{A}$ with $\bm{K}$, and solving for $\bm{y}_j$ at unknown indexes $j$}.

Finally, one needs to convert back from the encoded space to the original space using:
\begin{IEEEeqnarray}{c}
\underset{n\times1}{\bm{x}} = h\left(\underset{p\times1}{\bm{y}}\right)
\end{IEEEeqnarray}
where $h(\cdot)$ should approximate $g^{-1}(\cdot)$.

Note that the term \emph{extended dynamic mode decomposition (EDMD)} refers to this same technique, except it uses orthonormal polynomial basis functions as $g(\cdot)$. Koopman theory in more general and can work with any $g(\cdot)$. For example, one can use a neural network encoder-decoder architecture to represent $g(\cdot)$ and $h(\cdot)$.

\break

\bibliographystyle{plain}
\bibliography{koopman_theory}

\section*{Appendix: Variable reference}
$m$ -- Number of data measurements excluding the 0th measurement. Thus, total number of data measurements is $m+1$.

$n$ -- Number of original states.

$p$ -- Number of transformed or encoded states.

$r$ -- Rank.

\end{document}